\documentclass[twoside, 12pt]{article}
\usepackage[russian, english]{babel}
\usepackage{epsfig}
\usepackage{amssymb,amsmath,amsfonts,amsthm,enumerate}
\usepackage{amsfonts}
\usepackage{amsmath}
\usepackage{graphicx}
\usepackage{amsthm}
\usepackage[cp1251]{inputenc}
\usepackage[T2A]{fontenc}
\usepackage{mathrsfs}
\usepackage[english]{babel}
\newtheorem{theorem}{Theorem}[section]

\newtheorem{remark}{Remark}
\newtheorem{prop}{Proposition}[section]
\newtheorem{corollary}{Corollary}[section]
\numberwithin{equation}{section}
 \def\@evenhead{\vbox{\hbox to \textwidth{\thepage\hfil\sl\leftmark\strut}\hrule}}
 \def\@oddhead{\vbox{\hbox to \textwidth{\rightmark\hfill\thepage\strut}\hrule}}

\begin{document}
 \sloppy

\centerline{\bf Differentiation operator in the Beurling space of ultradifferentiable functions}
\centerline{\bf of normal type
 on an interval}     % Title of talk

\vskip 0.3cm

\centerline{\bf  N.~Abuzyarova}        % Authors from the same institution

\markboth{\hfill{\footnotesize\rm  N.~Abuzyarova }\hfill}
{\hfill{\footnotesize\sl  Differentiation operator in the Beurling space of ultradifferentiable functions
of normal type
 on an interval
}\hfill}
\vskip 0.3cm

\vskip 0.7 cm


\vskip 0.2cm

\vskip 0.2cm

\noindent {\bf Abstract.}
In this paper we study   closed 
subspaces of  ultradifferentiable functions which are invariant under the differentiation operator.
We propose a version of spectral synthesis which takes into account the presence of non-trivial differentiation invariant subspaces
containing no exponential monomials. As an application, we describe the sets of solutions of finite and infinite  systems of <<local>> homogeneous 
convolution equations.

\section{Introduction}

Let $\omega:[0;\infty) \to [0;\infty)$ be non-decreasing continuous  function satisfying the  conditions:

\noindent
the function $\varphi (t):=\omega (e^t)$ is convex on $[0;\infty)$,

\noindent 
$
\omega(x) =o (x),$ and
$\mathrm{ln}\, x =o(\omega (x))$  as $x\to\infty,$

\noindent
$\omega$ is an <<almost subadditive>> weight:
\begin{equation}
\forall \sigma>1 \ \exists C>0:\ \omega (x+y)\le\sigma (\omega(x)+\omega (y))+C, \ \forall x, y\ge 0 ,
\label{al-sa}
\end{equation}

\noindent
$\omega$ is {\it non-quasianalytic weight}:
\begin{equation}
\int_{1}^{\infty}\frac{\omega (x)}{x^2}\mathrm{d}x<\infty .
\label{non-qa}
\end{equation}

Let $\{ c_k\}$, $k=1,2,\dots ,$
be such that $0<c_k\nearrow a,$ $0<a\le+\infty .$ Then,
 $\{[-c_k;c_k]\}$ is the  
increasing sequence of closed intervals exhausting 
 $(-a;a)$.
Given $f\in C^{\infty} (-a;a),$ $q\in (0;1)$ and $k\in\mathbb N ,$ we set
\begin{equation}
\| f\|_{\omega, q,k} =
\sup_{j\in\mathbb N_0} \sup_{|x|\le c_k} \frac{|f^{(j)}(x)|}{e^{q\varphi^*(j/q)}},
\label{norma-udf}
\end{equation}
where $\varphi^*(s)=\sup\limits_{t} (st-\varphi (t)).$

{\it The space of ultradifferentiable functions (UDF) of normal type on } $(-a;a)$ is defined as
$$
\mathcal U_a =\{ f\in C^{\infty} (-a;a):\ \| f\|_{\omega, q,k}<\infty \  \forall q\in (0;1), \ \forall k=1,2,\dots \} .
$$
The system of semi-norms $\|\cdot\|_{\omega,q,k}$  defines  a locally convex topology in $\mathcal U_a.$
Equipped with this topology,  $\mathcal U_a$ becomes a locally convex space of  $(M^*)$-type.

For an arbitrary interval $(a;b)\subset\mathbb R,$ we denote by $\mathcal U (a;b)$ 
the space of UDF with shifted argument:  
$$f\in \mathcal U (a;b)\Longleftrightarrow 
f(x+(a+b)/2)\in \mathcal U_{(b-a)/2}.$$

Let $D=\frac{\mathrm{d}}{\mathrm{d}x}$ be the differentiation operator acting in
 $\mathcal U_a,$
$W\subset\mathcal U_a$ be a {\sl closed differentiation invariant} subspace,
 that is $D(W)\subset W$.
Shortly,   $W$ is said to be {\it $D$-invariant}
subspace.

It is well-known that  
{\sl root elements} of the differentiation operator are exponential monomials  
$t^ke^{-\mathrm{i}\lambda}t,$ $k\in\mathbb Z_+,$ $\lambda\in\mathbb C.$
We denote by  $\mathrm{Exp} (W)$ the set of all exponential monomials contained in  $W.$

Given an ultradistribution $S\in\mathcal U'_a$ 
 the  formula
$$g(y)=T_S(f)(y)=T_S(f(x+y))$$
defines 
a {\it convolution operator} $T_S: \mathcal U_a \to \mathcal U(c+a; a-d),$
		where  $[c;d] =\mathrm{ch}\, \mathrm{supp}\, S$ is 
	  convex hall of the support of $S.$
		The kernel of $ T_S$ provides a classical example of $D$-invariant subspace.

Let $S\in\mathcal U'_a.$
It is known that if 
its Fourier-Laplace transform $\mathcal F(S)$ is a divisor of  $\mathcal F(\mathcal U'_a)$
then every function
$f\in \mathrm{ker}\, T_S$
can be represented as series of exponential polynomials contained in $\mathrm{ker}\, T_S$.
The series  converges with grouping in $\mathcal U_a$.
This is the result  due to D.A. Abanina  \cite{Vl-AbDA}.
In the next sections, we  recall the explicit definitions of Fourier-Laplace transform
$\mathcal F$, multiplier and divisor of the space $\mathcal F(\mathcal U'_a)$.
At the moment, we only need to mention that  if $\mathcal F(S)$ is a divisor of  $\mathcal F(\mathcal U'_a)$
then $\mathrm{supp}\, S=\{ 0\}$ (see \cite{Vl-AbDA}).

On the other hand, for an ultradistribution  $S$ with an arbitrary compact support in $(-a;a)$,
we do not know any results even on  the approximation of functions in
$ W_S:=\mathrm{ker}\, T_S$ by the elements of $\mathrm{span}\, \mathrm{Exp} (W_S).$

Consider   more general situation.
Given  $A,\ B \subset\mathbb R$,  denote by $A\div B$  
their 
 {\it geometric difference}
consisting of all $x\in\mathbb R $ such that  $x+B\subset A.$

Let $S\in\mathcal U'_a$ and $I\subset (-a;a)$ be a relatively closed interval such that
$\widetilde{I}:=I\div \mathrm{ch}\, \mathrm{supp}\, S\neq \emptyset .$
We define {\sl <<local>>} convolution operator $T_{S,I}: \mathcal U_a\to \mathcal U (\widetilde{I})$  by 
  	setting
				$$
		g=T_{S,I} (f),
		\quad g(y)=S(f (x+y)),\ y\in \widetilde{I}.
		$$
		Clearly, the kernel $W_{S,I}:=\mathrm{ker}\, T_{S,I}$
		is a closed $D$-invariant 
		subspace in $\mathcal U_a.$
		It is also not difficult to see that
			$W_{S,I}$ contains all functions
 $f\in\mathcal U_a$ vanishing on 
 $I.$
 So, we face a  problem
because  any non-zero function $f\in\mathcal U_a$ vanishing on $I $ 
 cannot belong to the closure  of 
$\mathrm{span}\, \mathrm{Exp} (W_{S,I})$ in $\mathcal U_a.$
This is the consequence of two facts:

\noindent
 1) for any interval $X\subset\mathbb R$,
there are   {\sl continuous} embeddings
$$\mathcal U (X)\subset C^{\infty}(X)\subset C(X);$$

\noindent
2) for any relatively closed interval $X\subset (-a;a)$ and  mean periodic function $f$ on $X$,
 its mean periodic continuation to   $(-a;a)$ is unique
   (see \cite{Sedl}, \cite{Leont}).

Denote by $\mathcal U_{min}(\mathbb R)$ the space of UDF of minimal type on $\mathbb R$.
In  \cite{MTV}, the  authors 
consider a subspace $\mathcal W\subset \mathcal U_{min}(\mathbb R)$
which is similar to $W_{S,I}.$   
They 
prove 
that any function  in $\mathcal W$ can be {\sl locally} expanded in a series (with grouping)
of exponential polynomials contained in $\mathcal W.$ This series converges  in a small interval
under the condition that  
$\mathcal F (S)$  is a divisor 
of the space $\mathcal F (\mathcal U'_{min}(\mathbb R))$,
where $S$ is the ultradistribution defining $\mathcal W.$

In fact, the  above questions are 
particular cases 
of the following {\sl spectral synthesis problem}:
how to reconstruct $D$-invariant space $W\subset \mathcal U_a$ knowing  $\mathrm{Exp} (W)?$
Generally speaking, it is not equivalent to establish
the relation $W=\overline{\mathrm{span}\,\mathrm{Exp} (W)}.$

Given relatively closed interval $I\subset (-a;a)$,
we set 
$$
W_I=\{ f\in\mathcal U_a:\ f=0 \ \text{на}\ I \}.
$$
If  $I\neq (-a;a)$ then $W_I$ is non-trivial $D$-invariant subspace containing no
exponential monomials.
Following  \cite{Al-Kor}, we call  $W_I$ {\it residual subspace}. 
 We will show that 
any $D$-invariant subspace $W\subset\mathcal U_a$ 
has  
 {\it residual interval} $I_W$ (see Proposition
  \ref{pr-3} below). It is  defined by the following requirements:
 $I_W$ is relatively closed in $(-a;a)$, $W_{I_W}\subset W,$ and at the same time,
$W_{I}\setminus W \neq{\emptyset}$ for any relatively closed interval
   $I\subsetneq I_W$.  
  
Now,
taking into account the presence of residual subspaces, we formulate spectral synthesis problem  
for the differentiation operator $D$ in  $\mathcal U_a$
in a similar way as it was done for $D$ in $C^{\infty} (a;b)$
(see \cite{Al-Kor}, \cite{DAN-14}, \cite{ABB}, \cite{MZ}).

Let  $W\subset \mathcal U_a$ 
be a $D$-invariant subspace with  residual interval  
 $I_W$ and supplies of exponential monomials 
$\mathrm{Exp} (W) $. 
We say that $W$ {\it admits weak spectral synthesis} 
if
\begin{equation}
W=\overline{\mathrm{span} (\mathrm{Exp} (W))+W_{I_W}}.
\label{w-synt}
\end{equation}
  $D$-invariant subspace
 $W$ {\it admits spectral synthesis} if
\begin{equation}
W=\overline{\mathrm{span} (\mathrm{Exp} (W))}.
\label{synt}
\end{equation}
Clearly, in this case, $I_W =(-a;a).$
{\sl Spectral synthesis problem:} given $D$-invariant subspace $W\subset \mathcal U_a,$
to obtain the conditions under which the relation (\ref{w-synt}) holds?

In  \cite{DAN-14}, \cite{MZ}, 
we modified classical dual scheme to solve the similar 
 problem in $C^{\infty} (a;b)$.
This scheme goes back to L. Ehrenpreis and I.F. Krasichkov-Ternovskii \cite{Kr-Ter}.
There, we also partially used some results on spectral analysis and synthesis for 
the differentiation operator in $C^{\infty} (-a;a)$ due to A.Aleman and B. Koremblum \cite{Al-Kor}. These authors
got their results
by different methods and approaches.
 However, all of them may be  obtained 
by using our dual scheme.

In this paper, we study the spectral synthesis 
problem for the  differentiation operator in $\mathcal U_a$.
 The results  on the general theory of UDF and ultradistributions
 due to A.V. Abanin
(see \cite{Ab-m}, \cite{Ab-UMN}) make it possible to use  here the  dual approach too.
The solution of the spectral synthesis problem will be obtained for the most general situation.

We stick to the following  plan. 

In Section 2, we introduce a topological module 
of entire functions $\mathcal P_a$,
 describe duality between  $D$-invariant subspaces 
in $\mathcal U_a$
 and closed submodules in  $\mathcal P_a$,
establish  equivalence of the spectral synthesis problem in
 $\mathcal U_a$ and the dual problem of local description  in 
 $\mathcal P_a$.
In Section 3, we explore characteristics and properties of closed submodule
$\mathcal J\subset\mathcal P_a$  we will  need in. We  also solve the
local description problem there. 
It is equivalent to find the solution of the
 spectral synthesis  problem
 in $\mathcal U_a$.

\section{Duality of  $D$-invariant subspaces and weighted submodules of entire functions}

In this section, we establish duality between $D$-invariant subspaces of UDF and 
submodules of entire functions. It leads to the equivalence of the spectral synthesis problem for $D$-invariant
subspaces 
and the local description problem for
 submodules. 
We apply duality scheme mentioned in Introduction.
The classical version of this scheme deals with $D$-invariant subspaces of functions which are holomorphic
in a convex domain (see \cite{Kr-Ter}).
We have modified the scheme to study $D$-invariant subspaces in 
 $C^{\infty}(a;b)$ (see 
\cite{UMJ-14}, \cite{DAN-14}, \cite{MZ}). There, we have  also faced some difficulties 
because of presence of non-trivial
residual subspaces in $C^{\infty}(a;b)$.

According the general theory of UDF developed by A.V. Abanin ( \cite{Ab-m}, \cite{Ab-UMN}),
all polynomials and exponentials $e^{\mathrm{i}\lambda t},$ $\lambda\in\mathbb C,$ belong to 
$\mathcal U_a$; and this space forms a ring with respect the operation of multiplication of functions.
The space of ultradifferentiable test functions contains
cut-off functions and partition of unity.
The strong dual space  $\mathcal U'_a$ is a space of all ultradistributions
with compact supports in  $(-a;a).$
Both spaces, $\mathcal U_a$ and $\mathcal U'_a$ are reflexive because of  first of them is of $(M^*)$-type,
and the second is of  $(LN^*)$-type. 
It implies the following assertion.

\begin{prop}
({\sl General duality principle.})
Between the set of all closed subspaces  $\{ W\}$ of the space $\mathcal U_a$
and the set  of all closed subspaces $\{ V\}$ of the strong dual space $\mathcal U'_a$
there is one-to-one correspondence defined by the rule:
$$
W\longleftrightarrow V \Longleftrightarrow V=W^0,
$$
where  $W^0=\{ S\in \mathcal U'_a :\ S(f)=0 \ \forall f\in W \}.$
\label{pr-1}
\end{prop}

The subspace  $W^0$ is called {\it annihilator subspace} for $W$.

\smallskip

Let $c_k$ be the same as in (\ref{norma-udf}), $0<r_k\nearrow 1,$ . 
We set
$$
 P_k=\left\{ \varphi\in H (\mathbb C): \ \ \|\varphi \|_k=\sup\limits_{z\in\mathbb C} 
\frac{|\varphi (z)|}{e^{r_k\omega (|z|)+c_k|\mathrm{Im}\, z|}}<+\infty \right\},\quad k=1,2,\dots
$$
It is easy to check that  $P_k$ is a Banach space, $k =1,2, \dots$

Denote by  $\mathcal F$ the Fourier-Laplace acting in $\mathcal U'_a$:
$$
\mathcal F (S) (z)=S(e^{-\mathrm{i}tz}),\  z\in\mathbb C.
$$
It is known that $\mathcal F$
is a linear topological isomorphism between 
 $\mathcal U'_a$ and $\mathcal P_a:=\lim\, \mathrm{ind} P_k $ 
 (see \cite{Ab-m}, \cite{Ab-UMN}).

Because of (\ref{non-qa}),  any  $\varphi\in\mathcal P_a$ is
of the Cartwright class of entire functions.
In particular,  it  is an entire function of completely regular growth with respect to the order  1, 
and the indicator diagram of $\varphi$ equals $[\mathrm{i} h_{\varphi}(-\pi/2);\mathrm{i} h_{\varphi}(\pi/2)]$,
where $h_{\varphi} $ in the indicator function of $\varphi$.

By Fourier-Laplace transform, the generalized differentiation operator 
$D':\mathcal U'_a\to\mathcal U'_a$ 
transforms to the operator of multiplication by $(-\mathrm{i}z)$. The last one acts continuously in
$\mathcal P_a$.
It means that $\mathcal P_a$ is a topological module over the  ring $\mathbb C[z];$
and if $W\subset \mathcal U_a$ is a $D$-invariant subspace  then $\mathcal J=\mathcal F(W^0)$
 is a closed submodule in $\mathcal P_a.$

In our further considerations, 
we will always deal with  {\sl closed } submodules $\mathcal P_a$ and 
 will omit the word <<closed>> talking about them. 
\smallskip
Given $D$-invariant subspace $W$,
we denote by $\Lambda_W$ the sequence of multiple points 
 $\{(\lambda_k;m_k)\}$
such that 
 $$
\mathrm{Exp} (W) =\{ t^je^{-\mathrm{i}\lambda_k}t, \ j=0,1,\dots ,m_k-1\}_{k=1}^{\infty}.
$$
Denote by $\mathcal Z_{\varphi} $  zero set of  function $\varphi\in\mathcal P_a,$
and by $\mathcal Z_{\mathcal J}$ {\it  zero set  }  of a submodule $\mathcal J\subset\mathcal P_a,$
 that is $\mathcal Z_{\mathcal J}=\bigcap\limits_{\varphi\in \mathcal J}\mathcal Z_{\varphi}$.
In the other words, $(\mu_k;n_k)\in\mathcal Z_{\mathcal J}$ if and only if any function $\psi\in\mathcal J$
vanishes at $\mu_k$ with a  multiplicity at least $n_k$, and 
there exists $\varphi\in\mathcal J$ vanishing at $\mu_k$ with the  multiplicity $n_k.$

By Proposition \ref{pr-1}, taking  into account the isomorphism  $\mathcal P_a=\mathcal F(\mathcal U_a')$ and
 the relation between   differentiation  in  $\mathcal U_a$ and multiplication by
 $(-\mathrm{i}z)$ in $\mathcal P_a$, we get
 
\begin{prop}({\sl Special duality principle.})
Between the set of all $D$-invarinat subspaces  $\{ W\}$ 
and the set  of all submodules $\{ \mathcal J\}$ 
there is one-to-one correspondence defined by the rule:
$$
W\longleftrightarrow \mathcal J \Longleftrightarrow \mathcal J=\mathcal F(W^0);
$$
in addition, $\Lambda_W=\mathcal Z_{\mathcal J}.$
\label{pr-2}  
\end{prop}

The submodule $\mathcal J=\mathcal F(W^0)$ is said to be  an {\it annihilator submodule}
for $W$.

Given submodule $\mathcal J\subset \mathcal P_a$, 
 we define its {\it indicator segment} $[c_{\mathcal J};d_{\mathcal J}]\subset\overline{\mathbb R}$
 setting
$c_{\mathcal J} =\sup\limits_{\varphi\in\mathcal J} h_{\varphi}(-\pi /2),$  
$d_{\mathcal J} =\sup\limits_{\varphi\in\mathcal J} h_{\varphi}(\pi /2)$.

\begin{prop}
For any  $D$-invariant subspace $W\subset\mathcal U_a$, there exists its residual interval  $I_W$ 
equaled to $[c_{\mathcal J};d_{\mathcal J}] \bigcap(-a;a)$, where $\mathcal J=\mathcal F (W^0)$.
\label{pr-3}
\end{prop}

{\bf Proof.}

\smallskip

Set $I_0=(-a;a)\bigcap [c_{\mathcal J };d_{\mathcal J}].$

It is known that the shift operators  
$$
f \mapsto f(\cdot +y), \quad (f\mapsto f (\cdot-y)), \quad y>0,
$$
acts continuously  in $\mathcal U(-a;+\infty)$
and in $\mathcal U(-\infty ;a)$, respectively.
For  any $f\in W_{I_0}$, we  have the representation
$$
f=f_{-}+f_{+},\quad f_{-}\in W_{I_{-}},\ \ f_{+}\in W_{I_+},
$$
where  $I_{-}=(-\infty;d_{\mathcal J}],$ $I_+ =[c_{\mathcal J};+\infty ),$
$W_{I_{-}}\subset \mathcal U(-\infty ;a)$, $ W_{I_+}\subset\mathcal U (-a;+\infty).$
Also, for any $S\in\mathcal F^{-1} (\mathcal J),$
$$\mathrm{supp}\, g (\cdot-y)\bigcap \mathrm{supp}\, S=\emptyset,\quad
g\in W_{I_{-}}, \ y>0, $$
$$\mathrm{supp}\, \tilde{g} (\cdot+y)\bigcap \mathrm{supp}\, S=\emptyset,\quad
\forall\, \tilde{g}\in W_{I_{+}}, \ y>0. $$

From the above we derive  that
$$
S(f)=S(f_{-}+f_{+})=\lim_{y\to 0+} \left(S(f_{-} (x-y)) +S(f_+ (x+y))\right) 
=0
$$
for any $S\in\mathcal F^{-1} (\mathcal J)$.
The duality principle (Proposition \ref{pr-1}) implies that $W_{I_0}\subset W.$

Let an interval $I'\subsetneq I_0$ be  relatively closed in $(-a;a)$.
From the definitions of 
$c_{\mathcal J}$ and $d_{\mathcal J}$ and general theory of UDF and ultradistributions,
 we get that for any
 $c'\in (c_{\mathcal J };d_{\mathcal J})\setminus I',$
there exist
 $S\in\mathcal F^{-1}(\mathcal J),$  $f\in\mathcal U_a$
and $\delta >0$ satisfying
$$
S(f)\neq 0, \quad \mathrm{supp}\, f\subset (c'-\delta ;c'+\delta)\subset (c_{\mathcal J };d_{\mathcal J})\setminus I'.
$$
By the duality principle, we conclude that $f\not\in W.$ 
However, $f\in W_{I'}.$
It means that $I_0$ is the  minimal one among all intervals relatively closed in $(-a;a)$  intervals  with the property $W_{I}\subset W.$
That is,   $I_W=I_0.$

Q.E.D.

\medskip

We are going to show that the spectral synthesis
problem for a  $D$-invariant subspace $W\subset\mathcal U_a$ 
is equivalent to the problem of reconstructing of submodule $\mathcal J=\mathcal F (W^0)$ 
 under the assumption that its zero set and indicator segment are known.

Notice that  $D$-invariant subspace  $W$ admitting weak spectral synthesis (\ref{w-synt}),
 is minimal among all
  $D$-invariant subspaces $\widetilde{W}\subset\mathcal U_a$ such that
	$$
	I_{\widetilde{W}}=I_W,\quad \mathrm{Exp} (\widetilde{W})=\mathrm{Exp} (W).
	$$
	By Proposition  \ref{pr-2},  we see that the annihilator submodule $\mathcal J=\mathcal F (W^0)$
	is maximal among all submodules 
	$\widetilde{\mathcal J}\subset \mathcal P_a$ with
		$$
	\mathcal Z_{\widetilde{\mathcal J}}=\mathcal Z_{\mathcal J}\quad\text{and}\quad [c_{\widetilde{\mathcal J}};d_{\widetilde{\mathcal J}}] 
	=[c_{\mathcal J};d_{\mathcal J}].
	$$
		
		A submodule $\mathcal J\subset\mathcal P_a$
		is said to be  {\it weakly localisable}
		  if it has the described maximality
		property.
		Such submodule contains  {\sl all} functions $\psi\in\mathcal P_a$ 
		 with indicator diagrams 
	$ [\mathrm{i} c_{\psi}; \mathrm{i} d_{\psi}]\subset  [\mathrm{i} c_{\mathcal J};\mathrm{i} d_{\mathcal J}]$
	vanishing on  $\mathcal Z_{\mathcal J}$.
	
	We have proved the following assertion  (compare with \cite[Proposition 2]{MZ}).
	
		\begin{prop}
	$D$-invariant subspace $W\subset \mathcal U_a$ admits weak spectral synthesis if and only if its annihilator submodule
		$\mathcal J\subset \mathcal P_a$
	is weakly localisable.
	\label{pr-4}
\end{prop}

\begin{remark}
In particular, if  $I_W=(-a;a)$ then
 Proposition \ref{pr-4} establishes classical duality between the
spectral synthesis problem and the problem of local description of submodules.
 In this case, the submodule  $\mathcal J$ is called   {\sl localisable} or {\sl ample}
 (see \cite{Kr-Ter},
\cite{IF-loc-1}, \cite{DAN-14}, \cite{MZ}).
\end{remark}

\section{Local description of submodules in $\mathcal P_a$}

\subsection{Weak localisability of stable submodules.}

Let $\mathcal P$ be a locally convex space of entire functions.
It
is called {\it $b$-stable}
if, for any bounded subset $B\subset\mathcal P$,  
the set of all {\sl entire }  functions of the form 
$$ \psi =\frac{\varphi}{z-\lambda},\quad \lambda\in\mathbb C,  \quad \varphi\in B,$$
 is a bounded subset of $\mathcal P$.

As we have noticed above,    $\mathcal P_a$ is a locally convex space of $(LN^*)$-type.
It is well-known that a subset 
$B\subset \mathcal P_a$ is bounded if and only if $B$ is a bounded subset in
some $ P_k$ 
(see \cite[Theorem 2]{Seb}).
Taking into account this fact and the description of the topology  in $\mathcal P_a,$
we can easily derive that  $\mathcal P_a$ is
a bornological and  $b$-stable space. 
It means that we may apply
abstract methods developed 
in \cite{IF-loc-1},
\cite{IF-loc-2} to study  submodules in $\mathcal P_a.$

A submodule $\mathcal J\subset\mathcal P_a$  is {\it generated by functions} 
$\varphi_1,\dots, \varphi_m$ ({\it $m$-generated}) if  it is the closure of the set
$\{ p_1\varphi_1+\dots +p_m\varphi_m\},$ where $p_1,\dots, p_m$ are polynomials.

 A submodule $\mathcal J\subset\mathcal P_a$ 
 generated by one function
 $\varphi\in\mathcal P_a.$
 is said to be  {\it principal}.
 We denote it by $\mathcal J_{\varphi}.$

 Given a set $\mathcal S=\{\varphi_{\alpha}\}\subset\mathcal P_a$,
a submodule $\mathcal J\mathcal P_a$  {\it generated  by} $\mathcal S$
  equals the closure in $\mathcal P_a$
of the set of all finite sums
$$
p_1\varphi_{\alpha_1}+\dots +p_k\varphi_{\alpha_k},\quad
\text{where}\quad \varphi_{\alpha_j}\in\{\varphi_{\alpha}\},\quad p_{j}\in\mathbb C[z].
$$
That is, $\mathcal J$ is
minimal among all submodules containing $\mathcal S.$
 
\smallskip

Given a submodule $\mathcal J\subset\mathcal P_a,$
let $\lambda\in\mathbb C$ is a {\it  zero of $\mathcal J$ of multiplicity } $n_{\lambda}\ge 0$,
that is any function $\psi\in\mathcal J$ vanishes at $\lambda$ with the multiplicity at least $n_{\lambda}$
and there exists $\psi_{\lambda}\in\mathcal J$ for which $\lambda$ is a zero of multiplicity equaled to $n_{\lambda}.$
A submodule $\mathcal J$ is  {\it stable at the point } $\lambda$ 
if the following implication holds:
\noindent
if $\varphi\in\mathcal J$, and $\lambda$ 
is  a zero of $\varphi$ of the multiplicity $n>n_{\lambda}$
 then $\frac{\varphi}{z-\lambda}\in\mathcal J.$

A submodule $\mathcal J$ is  {\it stable } if it is stable at any point $\lambda\in\mathbb C.  $

From the results of  \cite[$\S 4$]{IF-loc-2} we derive that

\noindent
1) stability of a submodule $\mathcal J$
at one point $\lambda_0 \in\mathbb C$ implies its stability at any point $\lambda\in\mathbb C$;

\noindent
2) principal submodules in $\mathcal P_a$ are stable.

It is easy to see that any weakly localisable submodule
 $\mathcal J\subset\mathcal P_a$ is stable.
However, 

\noindent
a) there are stable submodules in $\mathcal P_a$ which are not weakly localisable (see remarks at the end of this subsection);

\noindent
b)  there are unstable submodules in $\mathcal J$ (Proposition \ref{pr-7} below).

We are going to obtain a condition guaranteeing weak localisability
of stable submodule.
In \cite{IF-loc-1}, I.F. Krasichkov-Ternovskii studied  abstract spaces of holomorphic vector functions.
For our purpose, we cite his results formulating them for  entire scalar functions.

Let  $ \mathcal P$ be a locally convex space of entire functions and a topological module over the ring  $\mathbb C[z].$ 
A closed submodule $ \mathcal J \subset \mathcal P$ is said to be {\it $b$-saturated with respect to a function} $\psi\in \mathcal P$
if there exists a bounded subset  
$B\subset\mathcal P$ such that the following implication holds: 

\noindent
$\rho $ is an entire function and
$$
 |\rho (z)\varphi (z)|\le |\varphi (z)|+|\psi (z)| ,\ \forall\, z\in\mathbb C ,
\quad \forall\, \varphi\in B\bigcap \mathcal J,
$$
  implies
 $\rho =const .$

Denote by $\mathcal P$  a bornological and $b$-stable space of entire functions.

{\bf Bornological version of the individual theorem.} (\cite{IF-loc-1})  
{\sl Let $\mathcal  J$ be a stable closed submodule in  $\mathcal P$ 
and a function $\psi\in  \mathcal P$  satisfy the condition
$\mathcal Z_{\mathcal J}\subset \mathcal Z_{\psi}.$
 
Then, $\psi\in \mathcal J$ if and only if  $\mathcal J$ is $b$-saturated with respect to $\psi.$}

\smallskip

Given $\varphi\in\mathcal P_a,$
we denote by
   ${\mathcal  J}(\varphi)$ the submodule
	consisting of all functions 
$\psi\in\mathcal P_a $ 
such that  $\psi =\rho\varphi,$ where  $\rho$ is an entire function of minimal type with respect to the order 1.
It means that the submodule
${\mathcal  J}(\varphi)$ contains all functions $\psi\in\mathcal P_a $ which zero sets satisfy
$\mathcal Z_{\varphi }\subset \mathcal Z_{\psi}$ and which  indicator diagrams equal the indicator diagram of $\varphi$.
Clearly, $\mathcal J(\varphi)$ is a weakly localisable submodule;
in particular, $\mathcal J_{\varphi}\subset \mathcal J(\varphi )$.  

 Now, we formulate and prove  Proposition \ref{pr-5}  and Theorem \ref{t-1}. They are very similar
to Lemma 1 and Theorem 1 in   \cite{MZ}.
There we considered the Schwartz module  $\mathbf P_a =\mathcal F ((C^{\infty} (-a;a))')$.
Proofs of both assertions are also similar to ones in \cite{MZ}.
Nevertheless, for the sake of  completeness, we state them here too.

\begin{prop}
Let  $\mathcal J\subset\mathcal P_a$ be a stable submodule.
Assume that  $\varphi\in\mathcal P_a$ satisfies the conditions:
$\mathcal Z_{\mathcal J}\subset\mathcal Z_{\varphi}$
and $[h_{\varphi}(-\pi/2); h_{\varphi}(\pi/2)]\subset  ( c_{\mathcal J};d_{\mathcal J})$. 

Then,  $J(\varphi)\subset \mathcal J.$
\label{pr-5}
\end{prop}  

\smallskip

{\bf Proof.}

Consider an arbitrary function   $\psi\in\mathcal J(\varphi).$
By the inequalities $c_{\mathcal J}<c_{\varphi},$ 
$d_{\mathcal J}>d_{\varphi}$
and the definitions of  $c_{\mathcal J}$ and $d_{\mathcal J}$, we get that
$$c_{\mathcal J}\le c_{\varphi_1}<c_{\varphi},\quad
d_{\varphi }< d_{\varphi_2}\le d_{\mathcal J} $$
 for some
 $\varphi_1,$
$\varphi_2\in\mathcal J$.

Set $\varphi_B =\varphi_1+\varphi_2.$ This is a function of completely regular growth with respect to the order 1.
 Notice that the indicator diagram of  $\psi\in\mathcal  J(\varphi)$ equals the closed interval
$\mathrm{i} [c_{\varphi};d_{\varphi}].$ Hence, it is compactly contained in the indicator
diagram   of $\varphi_B.$
It follows that
\begin{equation}
\frac{\psi (z)}{\varphi_B (z)}\to 0,\quad z=re^{i\theta}
\label{tend}
\end{equation}
as $r\to \infty$ staying outside of some set of zero respective measure.
The relation (\ref{tend}) holds uniformly with respect to $\theta \in \{ |\pi/2 -\theta|<\delta\}\bigcup \{|-\pi/2-\theta|<\delta\}$,
where $\delta>0$ is a sufficiently small fixed number.  

Show that the submodule  $\mathcal J$ is  $b$-saturated with respect to the function $\psi .$
We set $B=\{\varphi_B\}$ and
consider an entire function $\rho$ satisfying
 \begin{equation}
|\rho(z)\varphi_B(z)|\le |\psi (z) | +|\varphi_B (z)| ,\quad z\in \mathbb C.
\label{satur}
\end{equation}
Applying the theorem on summation of indicators 
 we derive 
that  $\rho$ is of minimal type with respect to the order 1.

By the
maximum modulus principle,
from (\ref{tend}) 
we get that
  $\rho$ is bounded on the imaginary axis.
It follows that $\rho=const ,$ and the stable submodule
$\mathcal J$ is $b$-saturated with respect to $\psi.$
By the bornoligical version of the individual theorem, we conclude that
$\psi\in \mathcal J$.

Q.E.D.

\medskip

\begin{theorem}
A stable submodule $\mathcal J\subset\mathcal P_a$ is weakly localisable if and  only if
it contains a function $\varphi$ with the property 
\begin{equation*}
\mathcal J(\varphi )\subset\mathcal J.
\end{equation*}
\label{t-1}
\end{theorem}

\smallskip

{\bf Proof.}

Clearly,  we only need to prove the sufficient part.

1) First, consider the situation, when $\mathcal J(\varphi)\subset \mathcal J,$ and the indicator diagram of $\varphi$ 
equals $\mathrm{i}\, [c_{\mathcal J};d_{\mathcal J}]$.
 (Unlike  in the Schwartz module \cite[Theorem 1]{MZ}, the case of $c_{\mathcal J}=d_{\mathcal J}$ is non-trivial in
$\mathcal P_a$.)

Let  $\psi\in \mathcal P_a$ be such that
  $$\mathcal Z_{\psi}\supset \mathcal Z_{\mathcal J},
\quad \text{and}  \mathrm{i}\,[c_{\psi };d_{\psi}]\subset
 \mathrm{i}\,[c_{\mathcal J};d_{\mathcal J}].$$

To make it sure that 
$\psi\in\mathcal J$, it is sufficient to check that 
 $\mathcal J$ is $b$-saturated with respect to $\psi.$ 

Set 
$$ 
B =\{ \widetilde{\varphi}\in H(\mathbb C): \ \ |\widetilde{\varphi}(z)|\le |\psi (z)|+|\varphi (z)|,\ \ z\in\mathbb C \}.
$$
By the topological properties of $\mathcal P_a$, 
 $B$ is  bounded.
Consider an entire function $\rho$  satisfying
 \begin{equation}
|\rho (z)\widetilde{\varphi}(z)|\le |\widetilde{\varphi} (z)|+|\psi (z)|, \ \ \forall \, z\in\mathbb C,
\quad \widetilde{\varphi}\in B\bigcap \mathcal J.
\label{satur1}
\end{equation}
In particular, we have
$$
|\rho (z)\varphi (z)|\le |\varphi (z)|+|\psi (z)|, \quad z\in\mathbb C.
$$
It follows that $\rho$ is of minimal type with respect to the order 1
and $\rho\varphi\in\mathcal J(\varphi)$.
  Taking into account the definition of $B$ and the relation $\mathcal J(\varphi)\subset \mathcal J, $
we obtain that
 $\rho\varphi \in B\bigcap \mathcal J.$
It means that (\ref{satur1}) holds with $\widetilde{\varphi} =\rho\varphi .$
 Hence,
$$
\left|\rho^2 (z) \varphi (z)\right|\le \left|\rho (z){\varphi} (z)\right|+|\psi (z)|\le 2 (|\varphi (z)|+|\psi (z)|) , 
\quad z\in\mathbb C.
$$
This estimate implies that 
$$
\frac{\rho^2}{2}\varphi\in B\bigcap \mathcal J.
$$
Now, setting  $\widetilde{\varphi} =\frac{\rho^2}{2}\varphi$ in (\ref{satur1}), 
we obtain
$$
\left|\frac{\rho^3 (z)}{2} \varphi (z)\right|\le  2 (|\varphi (z)|+|\psi (z)|) ,\quad z\in\mathbb C ,
$$
 and $$\frac{\rho^3}{2^2}\varphi \in B\bigcap \mathcal J.$$

Continuing to argue by the same way,
we prove that
$$
\frac{|\rho^n (z)|}{2^{n-1}} |\varphi(z)|\le |\varphi (z)|+|\psi (z)|,\quad z\in\mathbb C
$$
holds for all $n\in\mathbb N$.
Hence,  $\rho =const ,$
and   $\mathcal J$ is $b$-saturated with respect to $\psi .$

\smallskip

2) Now,  suppose that there exists $\varphi\in\mathcal P_a$ satisfying
 $$\mathcal J(\varphi)
\subset \mathcal J,\quad
 [c_{\varphi};d_{\varphi}]\subsetneq [c_{\mathcal J};d_{\mathcal J}]\subset (a;b).
$$ 
It means that at least one of the expressions, 
$\delta_1 =c_{\varphi}-c_{\mathcal J}$ or $\delta_2 =d_{\mathcal J}-d_{\varphi},$  is positive.

Consider in details the case when  both, $\delta_1$ and   $\delta_2,$ are positive.

By Proposition  \ref{pr-5}, 
 we have
 $$\mathcal J(e^{i\delta' z}\varphi )\subset \mathcal J, \quad
\mathcal J(e^{-i\delta'' z}\varphi) \subset\mathcal J$$
for all $\delta' \in [0;\delta_1)$
 and $\delta ''\in [0;\delta _2).$ 
 In particular,  
 \begin{equation}
 e^{i\delta ' z}\varphi , \ \ e^{-i\delta '' z}\varphi \in\mathcal J, \quad \delta'\in [0;\delta_1),\  \delta''\in [0;\delta_2). 
 \label{incl}
 \end{equation}

  Set $\Phi = (e^{i\delta_1 z}+e^{-i\delta_2 z})\varphi .$
 The relations
$$
\lim_{\delta'\to\delta_1}e^{i\delta ' z}\varphi =e^{i\delta_1 z}\varphi ,
\quad
\lim_{\delta''\to\delta_2}e^{-i\delta '' z}\varphi =e^{-i\delta_2 z}\varphi ,
$$ 
hold with respect to the topology of  $\mathcal P_a. $
 Together with  (\ref{incl}), they lead to the inclusion $\Phi\in\mathcal J.$
 
 Further, any function $\Psi\in\mathcal J(\Phi)$ is of the form
 $$\Psi=\rho \Phi =\rho (e^{i\delta_1 z}+e^{-i\delta_2 z})\varphi ,$$ 
 where $\rho$ is an entire function of minimal type with respect to the order 1.

One can easily check that  $\rho\varphi\in\mathcal P_a$. 
By this relation and Proposition  \ref{pr-5}, we derive that 
 $$\rho\varphi\in\mathcal J,
\quad
\Psi_{\delta'} =e^{i\delta ' z} \rho\varphi \in\mathcal J, \quad \forall\, \delta'\in (0;\delta_1),
 \quad \Psi_{\delta''}=e^{-i\delta '' z}\rho\varphi \in\mathcal J ,\quad  \forall\, \delta''\in (0;\delta_2) .
 $$
The relation  
$$\Psi=\lim\left(  \Psi_{\delta'}+\Psi_{\delta''} \right)\quad\text{as}\quad \delta'\to\delta_1,\ \ \delta''\to\delta_2,$$
leads to the inclusion $\Psi\in\mathcal J.$
Because of this is true for 
  an arbitrary function $\Psi\in\mathcal J (\Phi)$,
 we conclude that $\mathcal J(\Phi)\subset \mathcal J .$ 

The indicator diagram of the function $\Phi$ equals the indicator segment of 
 $\mathcal J.$
It means that we are in the conditions of the first part of the proof.

\medskip

3) It remains to consider the case when  $c_{\mathcal J}=a$ or (and) $d_{\mathcal J}=b$.

Let  $\Psi\in\mathcal P_a $ and
 $\mathrm{i} [ c_{\Psi};d_{\Psi}]\subset\mathrm{i} [c_{\mathcal J};d_{\mathcal J}]$, 
$\mathcal Z_{\Psi}\supset\mathcal Z_{\mathcal J} .$
Prove that $\Psi\in\mathcal J.$

For this purpose, we choose and fix a closed interval $[c';d']$ such that
\begin{equation}
[c';d']\subset (a;b)\bigcap [c_{\mathcal J};d_{\mathcal J}],\quad
[c_{\Psi};d_{\Psi}]\subset [c';d'],\quad [c_{\varphi};d_{\varphi}]\subset [c';d'] .
\label{ind-segm-2}
\end{equation}
Denote by $\mathcal J'$ a weakly localisable submodule with the indicator segment equaled to
$[c';d']$ and zero set $\mathcal Z_{\mathcal J'}=\mathcal Z_{\mathcal J}.$
It is easy to see that 
$\widetilde{\mathcal J} =\mathcal J\bigcap\mathcal J'$ 
is the  closed stable submodule
with the indicator segment  $[c';d']$ and zero set
 $\mathcal Z_{\widetilde{\mathcal J}}
=\mathcal Z_{\mathcal J}.$

From (\ref{ind-segm-2}), it follows that $\mathcal J (\varphi)\subset\widetilde{\mathcal J}.$
Applying first two parts of the proof, we conclude that $\widetilde{\mathcal J} =\mathcal J'.$ 
Again by   (\ref{ind-segm-2}), we obtain the required relation 
$$\Psi\in\widetilde{\mathcal J}\subset \mathcal J.$$

Q.E.D.

\medskip

\begin{corollary}
Let 
$\mathcal J\subset \mathcal P_a$ 
be a stable submodule.
If there exists a function $\varphi_0\in\mathcal J_0$ such that
the principal submodule generated by $\varphi_0$ is weakly localisable
then $\mathcal J$ is also weakly localisable.
\label{cor-1}
\end{corollary}

{\bf Proof.}

Clearly, we have
$$\mathcal J(\varphi_0)=\mathcal J_{\varphi_0}\subset\mathcal J.$$ 
By  Theorem  \ref{t-1},  $\mathcal J$ is weakly localisable.

Q.E.D.

\medskip

Given sequence  $\Lambda=\{(\lambda_j;m_j)\}$
of multiple points, its
{\it the completeness radius} $r(\Lambda )$
  is defined to be the infimum of the set of all positive numbers $r$ such that
the system $\{ t^k e^{-\mathrm{i}\lambda_j t}, \ k=0,\dots , m_j-1\}_{j=1}^{\infty},$
is  incomplete in $C^{\infty} (-r;r)$ (equivalently, in $C(-r;r)$).

There are a few equivalent geometric characteristics of $\Lambda.$
 First of them is Beurling-Malliavin density
$D_{BM}(\Lambda )$ (see \cite{IF-BM}, \cite{Koosis-2}).
It is well-known that $r(\Lambda)=\pi D_{BM} (\Lambda)$ ({\sl Beurling-Malliavin theorem on the radius of completeness}).
Taking into account this fact and another theorem due to Beurling  and
Malliavin ({\sl Multiplier theorem}, see \cite[X-XI]{Koosis-2}), we obtain

\begin{prop}
If a submodule $\mathcal J\subset \mathcal P_a $  has the property
 $d_{\mathcal J}-c_{\mathcal J} <2\pi D_{BM} (\mathcal Z_{\mathcal J})$
then $\mathcal J=\{ 0\} .$ 
\label{pr-6}
\end{prop}

\begin{theorem}
Let $\mathcal J\subset P_a $ be a stable submodule.

\noindent
If $d_{\mathcal J}-c_{\mathcal J} >2\pi D_{BM} (\mathcal Z_{\mathcal J}),$
then $\mathcal J$ is non-trivial  and weakly localisable.
\label{t-2}
\end{theorem}
 
{\bf Proof.}

By the above cited results due  to Beurling and Malliavin and Paley-Wiener-Schwartz theorem
\cite[теорема 7.3.1]{Horm}, taking into account the conditions on the weight $\omega,$
we derive that there exists
non-zero function $\varphi_0\in\mathcal P_a$ with the properties
$\mathcal Z_{\mathcal J}\subset \mathcal J_{\mathcal \varphi_0}$
and $[h_{\varphi_0}(-\pi/2);h_{\varphi_0} (\pi/2)]\subset (c_{\mathcal J};d_{\mathcal J}).$

From Proposition \ref{pr-5}, it follows that $\mathcal J(\varphi_0)\subset\mathcal J.$ 
 Now, Theorem  \ref{t-1} implies the required assertion.

Q.E.D.

\medskip

It is of interest that
the similar result is true for stable  submodules in the Schwartz module
$\mathbf P_a$ (see \cite[Theorem 2]{MZ}).

\begin{corollary}
If the indicator segment of stable submodule $\mathcal J\subset\mathcal P_a$ is  non-compact in $(-a;a)$
then $\mathcal J$ is weakly localisable.

In particular, a stable submodule $\mathcal J\subset\mathcal P_a$ is localisable
if and only if
$$c_{\mathcal J}=-a, \quad d_{\mathcal J}=a.$$
\label{cor-2}
\end{corollary}

Recall that an entire function $\psi_0$ is  a {\it  multiplier} of the space $\mathcal P_a$
if the correspondence $\varphi\mapsto\varphi\psi_0 $ defines a continuous  map 
 of $\mathcal P_a$ into itself.
 
According \cite[Proposition 1]{AbDA-12}, the set $\mathcal M_a$ of all multipliers of $\mathcal P_a$
has the  description:
\begin{equation}
\mathcal M_a=\left\{ \psi\in H (\mathbb C) :\ \ \forall\varepsilon>0 \ \ \sup_{z\in\mathbb C}
\frac{|\psi (z)|}{e^{\varepsilon \omega (z)+\varepsilon |\mathrm{Im}\, z|}}<\infty .
\right\}
\label{descr-mult}
\end{equation}

  A {\it divisor} of  $\mathcal P_a$ is a function $\psi_0\in\mathcal M_a$ such that the implication
$$
\Psi\in \mathcal P_a, \quad \frac{\Psi}{\psi_0}\in H(\mathbb C)\Longrightarrow \frac{\Psi}{\psi_0}\in\mathcal P_a
$$
holds.

\begin{remark}
	Because of Theorem \ref{t-2}, Corollary \ref{cor-1} is worth considering only for non-trivial
	submodule $\mathcal J$ satisfying the relation 
	\begin{equation}
		d_{\mathcal J}-c_{\mathcal J} =
		2\pi D_{BM} (\mathcal Z_{\mathcal J}).
		\label{crit}
	\end{equation}
Principal submodules generated by slowly decreasing functions 
(see \cite{UMJ-2016-1})
or by divisors of the space  $\mathcal P_a$ (see \cite{AbDA-12})
 satisfy (\ref{crit}), and they are weakly localisable.
It may be proved by the scheme
we applied in
\cite{NF-VINITI} together with the  
results on weighted polynomial approximation
obtained in  \cite[VI.H.2]{Koosis-1}.
\label{rem-1}
\end{remark}

 It is known  that there exist stable submodules satisfying (\ref{crit}) in the Schwartz module $\mathbf P_a$
which  are not  weakly localisable
(see \cite{MZ}, \cite{ABB}, \cite{UMJ-2016-1}).
We guess that the same is true in $\mathcal P_a.$
We are going to study in details the critical case (\ref{crit}) for submodules in $\mathcal P_{a}$
later.

%%%%%%%%%%%%%%%%%%%%%%%%%%%%%%%%%%%%%%%%%%%%%%%%%%%%%%%%%%%%%%%%%%%%%%%%%%%%%%%%
%%%%%%%%%%%%%%%%%%%%%%%%%%%%%%%%%%%%%%%%%%%%%%%%%%%%%%%%%%%%%%%%%%%%%%%%%%%%%%%

\subsection{Stability}

Theorems \ref{t-1} and \ref{t-2} and their corollaries motivate to
study conditions of stability for 
submodules in
$\mathcal P_a$.
As we have noticed at the beginnig of the previous subsection,
any principal submodule is stable.
 However, for 2-generated submodule, it may fail.

\begin{prop}
Let $\mathcal J$ be a submodule generated in $\mathcal P_a$ by functions $\varphi_1,$ $\varphi_2\in\mathcal P_a,$ 
and $\mathrm{i}[c_1; d_1],$
$\mathrm{i}[c_2; d_2]$ be
their indicator diagrams.

If   
   $[c_1; d_1] \subset (-a;0),$
$[c_2; d_2] \subset (0;a)$
then the submodule  $\mathcal J$ 
is not stable.
\label{pr-7}
\end{prop} 

{\bf Proof.}

By Proposition  \ref{pr-2},
there exists a unique  $D$-invariant subspace
$W\subset\mathcal U_a$ such that 
 $\mathcal F (W^0) =\mathcal J.$
Its residual interval equals  $[c_1;d_2]$ (Proposition \ref{pr-3}).

Given a function $f\in\mathcal U_a$, by the duality principles (Propositions \ref{pr-1}, \ref{pr-2}),
we see that
 $$f\in W 
\Longleftrightarrow S_j (f^{(k)})=0,\ \ k=0,1,2,\dots, \ \ j=1,2,$$
where  
 $S_j=\mathcal F^{-1}(\varphi_j)$
$j=1,2 .$ 
In particular, for an arbitrary  $\varepsilon>0$ if  
$f\in \mathcal U_a$ vanishes on the set 
 $[c_1-\varepsilon;d_1+\varepsilon]\bigcup [c_2-\varepsilon;d_2+\varepsilon]$
then $f\in W.$

Notice that $\mathcal Z_{\mathcal J}=\mathcal Z_{\varphi_1}\bigcap \mathcal Z_{\varphi_2},$
 and $2\pi D_{BM} (\mathcal Z_{\mathcal J})<d_2-c_1.$ 
Hence, if $\mathcal J$ is a stable submodule then  
it is weakly localisable (by Proposition \ref{pr-4} and  Theorem \ref{t-2}),
that is
$$
W=\overline{W_{I_W}+\mathrm{span} (\mathrm{Exp} (W)}.
$$
The set of exponents of the  system $\mathrm{Exp} (W)$ equals $\mathcal Z_{\mathcal J}.$
This system is incomplete in  $C[c_j;d_j],$ $j=1,2.$
Taking into account these two facts, consider a  function
  $g\in W$ vanishing on
	$[c_1;d_1]\bigcup [c_2;d_2]$,
but $g\neq 0$ on
	 $I_W$. 
	By	the above, $g$ belongs to the closure of $\mathrm{span} (\mathrm{Exp} (W)$
in each  of the spaces $C[c_1;d_1], $ $C[c_2;d_2]$, $C[c_1;d_2].$ 
 It follows that
 $g=0$  $I_W=[c_1;d_2]$ 
according to the result on  uniqueness of mean periodic continuation
(\cite[$\S\, 1$]{Sedl}, 
\cite[$\S\, 9$]{Leont}).

It leads to the contradiction. Hence,
the submodule  $\mathcal J$ is not stable.

Q.E.D.

\medskip

Below, we consider the situations, when finitely generated submodule is necessarily stable.

\begin{theorem}
Let $\varphi_1,$ $\varphi_2\in\mathcal P_a$  be functions of minimal type with respect to
the order 1, and $\varphi_1\varphi_2\in\mathcal P_a.$

Then, they  generate  stable submodule $\mathcal J \subset\mathcal P_a$.
\label{t-3}
\end{theorem}

{\bf Proof.}

Without loss of generality ({\sl wlog}), we may assume that 
$\varphi_1(0)=\varphi_2(0)=1.$

It is known that
$\mathcal J$ is stable if and only if the identically zero function is contained in the closure in $\mathcal P_a$
of the set $\{p\varphi_1-q\varphi_2\},$ where $p,$ $q\in\mathbb C[z],$
 and $p(0)=q(0)=1$ (\cite[Proposition 4.9]{IF-loc-2}).

Because of the condition $\varphi_1\varphi_2\in\mathcal P_a,$
there exist $r_1,$ $r_2\in (0;1)$ such that
$r_1+r_2<1$ and for any $\varepsilon>0$
$$
\sup_{z\in\mathbb C}\frac{\varphi_j(z)}{e^{r_j\omega (|z|)+\varepsilon|\mathrm{Im}\, z|}}<+\infty,\quad j=1,2.
$$
Fix an arbitrary  $\delta>0$ satisfying  $r_1+r_2+2\delta <1.$
Taking into account the initial properties of  $\omega$, 
we apply the theorem on weighted polynomial approximation on the real line for the restrictions of  entire functions of minimal type
 \cite[VI.H.2]{Koosis-1}. 
By this result, we find two sequences of polynomials, $\{p_k\},$
$\{q_k\}, $  such that $p_k(0)=1,$ $q_k (0)=1,$ $k=1,2,\dots $,
and
\begin{eqnarray}
\sup_{x\in\mathbb R}\frac{|\varphi_1(x)-p_k(x)|}{e^{(r_1+\delta)\omega (|x|)}}\to 0, \ \ k\to\infty,
\label{pol-appr-1}
\\
\sup_{x\in\mathbb R}\frac{|\varphi_2(x)-q_k(x)|}{e^{(r_2+\delta)\omega (|x|)}}\to 0, \ \ k\to\infty .
\label{pol-appr-2}
\end{eqnarray}
Now, by standart argument including
the application of Phragmen-Lindel\"of principle (see \cite[VI.H.2]{Koosis-1})
and the  criterion of convergency for a countable
sequence in locally convex space of  $(LN^*)$-type (see \cite{Seb}),
 we derive that 
$$q_k\varphi_1-p_k\varphi_2\to 0,\ \ k\to\infty,
$$
 with respect to the topology of
$\mathcal P_a.$
By the above cited criterion of stability, 
we conclude that
 $\mathcal J$ is stable submodule.

Q.E.D.
 
\medskip

\begin{corollary}
Let $\varphi_1,\dots, \varphi_m\in\mathcal P_a$ 
be functions of minimal type with respect to the order 1 and
\begin{equation}
\varphi_j\varphi_{j+1}\in\mathcal P_a,\ \ j=1,\dots, m-1.
\label{in-1}
\end{equation}

 Then, 

\noindent 
1) the submodule $\mathcal J$ generated by 
 $\varphi_1,\dots, \varphi_m$ in $\mathcal P_a$
is stable;

\noindent
2) if, in addition, 
there is a divisor of $\mathcal P_a$ among functions $\varphi_j,$ $j=1,\dots , m,$ 
then the submodule $\mathcal J$ is weakly localisable.
\label{cor-3}
\end{corollary}

{\bf Proof.}

1)  {\sl Wlog},
$$\varphi_j(0)=1,\ \ 
j=1,2,\dots, m.$$

From (\ref{in-1}) and Theorem \ref{t-3}, it follows that
2-generated submodules
 $$\mathcal J_{\varphi_j,\varphi_{j+1}}=\overline{\{p\varphi_j+q\varphi_{j+1}\}},\quad
p,\ q\in\mathbb C[z],\ \ 
 j=1,\dots , m-1,$$
are stable.
Applying the necessary part of the stability criterion
 (\cite[Proposition 4.9]{IF-loc-2}) for these submodules,
 we find generalizes sequences of polynomials
$p^{(j)}_{\alpha}$, $q^{(j)}_{\alpha}$ 
such that  
\begin{equation}
p^{(j)}_{\alpha}(0)=q^{(j)}_{\alpha}(0)=1,
\quad p^{(j)}_{\alpha}\varphi_j-q^{(j)}_{\alpha}\varphi_{j+1}\to 0, \quad \alpha\nearrow
\label{str}
\end{equation}
 
By the sufficient part of the same criterion, we see
that 
 $\mathcal J$ is stable if
for any complex numbers $c_j,$
$j=1,\dots ,m, $
satisfying  $c_1+\dots +c_m=0,$
identically zero function belongs to the closure of the set
$\{s_1\varphi_1+\dots +s_m\varphi_m\},$
where
\begin{equation*}
s_j\in\mathbb C[z],  \ \
s_j (0)=c_j,\quad
j=1,\dots , m.
%\label{usl-s-j}
\end{equation*}

Set
\begin{eqnarray*}
s_{1,\alpha}=c_1p^{(1)}_{\alpha},\quad
s_{2,\alpha}=(c_1+c_2)p^{(2)}_{\alpha}-c_1q^{(1)}_{\alpha},\\
s_{3,\alpha}=(c_1+c_2+c_3)p^{(3)}_{\alpha}-(c_1+c_2)q^{(2)}_{\alpha},\\
\dots, \\
s_{m-1,\alpha}=(c_1+\dots +c_{m-1})p^{(m-1)}_{\alpha}-(c_1+\dots +c_{m-2})q^{(m-2)}_{\alpha},\\
s_{m,\alpha}=c_m q^{(m-1)}_{\alpha} .
\end{eqnarray*}
It is easy to check that
$$
s_{j,\alpha} (0)=c_j, \ \forall\, \alpha,\ \ j=1,\dots , m,
$$
and
$$
s_{1,\alpha}\varphi_1+\dots +s_{m,\alpha}\varphi_m\to 0
$$
with respect to the topology of $\mathcal P_a$.
It follows that  $\mathcal J$ is a stable submodule.

\smallskip

2) Assume that $\varphi_1$ is a divisor of $\mathcal P_a.$
According to Remark \ref{rem-1}, the principal submodule $\mathcal J_{\varphi_1} $ is weakly localisable,
that is $\mathcal J_{\varphi_1}=\mathcal J (\varphi_1)\subset\mathcal J.$
By the first assertion of the corollary,
 $\mathcal J$ is a stable submodule. Applying  Corollary  \ref{cor-1},
we conclude that $\mathcal J$ is a weakly localisable submodule.

Q.E.D.

\medskip

\begin{corollary}
Let $\{\varphi_{\alpha}\}\subset\mathcal M_a$.

Then,

\noindent
1) $\{\varphi_{\alpha}\}$ 
generates a stable submodule $\mathcal J$ in
 $\mathcal P_a;$
 
\noindent
2) if, in addition,  $\{\varphi_{\alpha}\}  $ contains a divisor of $\mathcal P_a$
then  $\mathcal J$  is weakly localisable.
\label{cor-5}
\end{corollary}

{\bf Proof.}

{\sl Wlog,} we assume that $\varphi_{\alpha} (0)=1$ for any $\alpha.$ 
Let   $\psi\in\mathcal J$ and
 $\psi (0)=0.$

Every function in $\mathcal J$
 is of minimal type with respect to the order  1. 
Taking into account the description of $\mathcal M_a$ cited above 
and  Theorem \ref{t-3},
we conclude that the functions $\psi$ and $\varphi_{\alpha}$ generate 
stable submodule $\mathcal J_{\psi,\varphi_{\alpha}}$.
It leads to the relations
$$\frac{\psi}{z}\in\mathcal J_{\psi,\varphi_{\alpha}}\subset
\mathcal J.$$

2) The second assertion is proved by the similar way as it 
was done to justify the second part of Corollary \ref{cor-3}. 

Q.E.D.

\subsection{Applications to  $D$-invariant subspaces}

Now, we apply the previous results to study weak spectral synthesis
for the differentiation operator in $\mathcal U_a.$

By Propositions  \ref{pr-3}, \ref{pr-4}, Corollary \ref{cor-1} and Theorem \ref{t-2}, 
we obtain

\begin{theorem}
Let  $W\subset\mathcal U_a$ be $D$-invariant subspace with residual interval $I_W$
of length  $d$, and  $\Lambda$ be a sequence of exponents of
  $\mathrm{Exp} (W)$.

Assume that the annihilator submodule $\mathcal J=\mathcal F (W^0)$ is stable, and one of the
 conditions 

\noindent
a) $2\pi D_{BM } (\Lambda)<d$ 

or 

\noindent
b)  there exists $\varphi_0\in\mathcal J$ generating   weakly localisable
 principal submodule 
$\mathcal J_{\varphi_0}$

\noindent
hold.

Then, $W$ admits weak spectral synthesis.
\label{t-4}
\end{theorem}

\smallskip

\begin{corollary}
Let $S\in\mathcal U'_a$, and $I\subset (-a;a)$ 
is a relatively closed interval such that
 $I\div \mathrm{ch}\, \mathrm{supp}\, S\neq\emptyset.$

Then, $W=\mathrm{ker}\, T_{S,I}$ 
is $D$-invariant
subspace admitting weak spectral synthesis
with
$I_W=I;$ and
the sequence of exponents of 
$\mathrm{Exp}(W)$ equals
$\mathcal Z_{\varphi},$  where $\varphi=\mathcal F(S).$
\label{cor-6}
\end{corollary}

{\bf Proof.}

By standard methods, one can easily check that 
$W=\mathrm{ker}\, T_{S,I}$ is a $D$-invariant subspace.
Applying the duality principles 
 (Propositions \ref{pr-1}, \ref{pr-2}), 
 we justify that 
\begin{equation}
\mathcal J=\mathcal F(W^0) =\overline{\mathrm{span} \, 
\{ e^{-\mathrm{i}tz}p\varphi:\ t\in I\div\mathrm{ch}\,\mathrm{supp}\, S, \ p\in\mathbb C[z]\} },
\label{ann-sub}
\end{equation}
is the annihilator submodule of $W.$
It follows  that the indicator segment of
  $\mathcal J$ is $\overline{I},$
and its zero set  equals   $\mathcal Z_{\varphi}.$
Together with Propositions  \ref{pr-2}, \ref{pr-3}, these facts imply that
$I_W=I,$ and
the sequence of exponents of 
$\mathrm{Exp}(W)$ is $\mathcal Z_{\varphi}.$ 
 
Remembering that
$I\div \mathrm{ch}\, \mathrm{supp}\, S\neq\emptyset,$ 
we apply 
Beurling-Malliavin theorem on radius of completeness and Paley-Wiener-Schwartz theorem for ultradistributions.
It gives  us the inequality
$2\pi D_{BM } (\mathcal Z_{\mathcal J})<d$, where   $d$ is the length of $I$.
Now, by Theorem
 \ref{t-4}, 
 to justify that $W$ admits weak spectral synthesis it is sufficient 
to check that  $\mathcal J$ is a stable submodule.

{\sl Wlog}, we assume that $\varphi (0)=1$ and prove that $\mathcal J$ is stable at the origin.
Together with \cite[Proposition 4.2]{IF-loc-2} 
it will lead to the required assertion.

Consider $\Psi\in\mathcal J,$ $\Psi (0)=0.$
By (\ref{ann-sub}), we have
$$
\Psi =\lim_{\alpha\nearrow}\Psi_{\alpha}, 
$$
where
$$
 \Psi_{\alpha} =(a_{1,\alpha}e^{-\mathrm{i}t_{1,\alpha}z}p_{1,\alpha}+\dots
+a_{n_{\alpha},\alpha}e^{-\mathrm{i}t_{n,\alpha}z}p_{n_{\alpha},\alpha})\varphi,
$$
$ a_{j,\alpha}\in\mathbb C,$ $ t_{j,\alpha}\in I\div\mathrm{ch}\,\mathrm{supp}\, S,$
 $ p_{j,\alpha}\in\mathbb C[z] ,$
$ j=1,\dots, n_{\alpha}.$

Notice that
\begin{equation}
\Psi_{\alpha}(0)=(a_{1,\alpha}p_{1,\alpha}(0)+\dots
+a_{n_{\alpha},\alpha}p_{n_{\alpha},\alpha} (0))\to 0.
\label{h1}
\end{equation}
Setting
$$
\Phi_{\alpha} =\sum_{j=1}^{n_{\alpha}}a_{j,\alpha}e^{-\mathrm{i}t_{j,\alpha}z}
(p_{j,\alpha} -p_{j,\alpha} (0))\varphi +
\sum_{j=1}^{n_{\alpha}}a_{j,\alpha}p_{j,\alpha} (0)
(e^{-\mathrm{i}t_{j,\alpha}z}-1)\varphi,
$$
we get the representation
$$
\Psi_{\alpha}=\Phi_{\alpha} +\Psi_{\alpha} (0)\varphi .
$$ 
Clearly, $\Phi_{\alpha} (0)=0$ and
$$
\frac{\Phi_{\alpha}}{z}=
\sum_{j=1}^{n_{\alpha}}a_{j,\alpha}e^{-\mathrm{i}t_{j,\alpha}z}
q_{j,\alpha}\varphi +
\sum_{j=1}^{n_{\alpha}}a_{j,\alpha}p_{j,\alpha} (0)
\frac{e^{-\mathrm{i}t_{j,\alpha}z}-1}{z}\varphi,
$$
where $q_{j,\alpha}=\frac{p_{j,\alpha} -p_{j,\alpha} (0)}{z}$ are polynomials.

Further,  by the results of \cite[VI.H.1]{Koosis-1}, we have
$$
\frac{e^{-\mathrm{i}t_{j,\alpha}z}-1}{z}\varphi\in
\overline{\mathrm{span}\, \{e^{-\mathrm{i}tz}\varphi \ \ t\in I\div\mathrm{ch}\,\mathrm{supp}\, S\}},
$$
$j=1,\dots, n_{\alpha}$, $\forall\,\alpha,$
where the closure is taken with respect to the topology of $\mathcal P_a$.
Hence, $\frac{\Phi_{\alpha}}{z}\in \mathcal J$ for every $\alpha.$ 
 
Because of (\ref{h1}), we have
$\Psi=\lim\Phi_{\alpha}.$

From the  topological properties of $\mathcal P_a$ it is not difficult to
derive that
$$
\frac{\Psi}{z}=\lim \frac{\Phi_{\alpha}}{z},
$$
It follows that $\frac{\Psi}{z}\in\mathcal J$,
and
 $\mathcal J$ is a stable submodule.

Q.E.D.

\medskip

\begin{corollary}
Let 
$I\subset (-a;a)$ 
be a relatively closed interval,
\begin{equation}
S_j\in\mathcal U'_a,\quad  \mathrm{supp}\, S_j=\{ 0\},\quad
j=1,\dots, m,
\label{starr}
\end{equation}
 and
$S_j*S_{j+1}\in\mathcal U_a'$
$j=1,\dots, m-1.$

Then, $W=\bigcap\limits_{j=1}^{m} \mathrm{ker}\, T_{S_j, I}$ 
admits weak spectral synthesis,
$I_W=I,$ and 
the sequence of exponents of 
$\mathrm{Exp}(W)$ is $\bigcap\limits_{j=1}^{m}\mathcal Z_{\varphi_j},$  where $\varphi_j=\mathcal F(S_j).$
\label{cor-7}
\end{corollary}

{\bf Proof.}

{\sl Wlog,} we assume that $0\in I$ and
 $\varphi_j(0)=1$, $j=1,\dots, m.$

By standard methods including the duality principles
(Propositions \ref{pr-1}, \ref{pr-2}), we  derive  that
 $W$ is a
$D$-invariant subspace in $\mathcal U_a$ and
\begin{equation}
\mathcal J=\mathcal F(W^0) =\overline{\mathrm{span} \, 
\{ e^{-\mathrm{i}t_jz}p_j\varphi_j:\ t_j\in I, \ p_j\in\mathbb C[z], \ j=1,\dots, m\} }
\label{ann-sub-2}
\end{equation}
is its annihilator submodule.

From (\ref{ann-sub-2}), it follows that $\overline{I}$ is the indicator segment of  $\mathcal J$,
and $$\mathcal Z_{\mathcal J} =\bigcap\limits_{j=1}^{m}\mathcal Z_{\varphi_j}.$$
It imples that $I_W=I,$ and the sequence of exponents of $\mathrm{Exp}(W)$ 
equals 
$\mathcal Z_{\mathcal J} =\bigcap\limits_{j=1}^{m}\mathcal Z_{\varphi_j}.$ 

By 
Beurling-Malliavin theorem on radius of completeness, Paley-Wiener-Schwartz theorem for ultradistributions
and  the condition (\ref{starr})
we obtain that
$2\pi D_{BM } (\mathcal Z_{\mathcal J})<d$, where  $d$ is the length of $I$.

Because of Theorem \ref{t-4} and \cite[Proposition 4.2]{IF-loc-2},
 the only thing we need to do is to check that
$\mathcal J$ is stable at the origin.

Let $\Psi\in\mathcal J$ and $\Psi (0)=0.$
From  (\ref{ann-sub-2}) it follows that
$$
\Psi =\lim\Psi_{\alpha},
$$
where  
\begin{multline*}
\Psi_{\alpha} =\sum_{j=1}^{m} a_{j,\alpha}e^{-\mathrm{i}t_{j,\alpha}z} (p_{j,\alpha}-p_{j,\alpha} (0))\varphi_j
+\sum_{j=1}^{m} a_{j,\alpha}p_{j,\alpha} (0)(e^{-\mathrm{i}t_{j,\alpha}z}-1) \varphi_j+\\
+\sum_{j=1}^{m} a_{j,\alpha}p_{j,\alpha} (0)( \varphi_j-\varphi_1) +\psi_{\alpha}.
\end{multline*}
Notice that
$$
\psi_{\alpha} =\left(\sum_{j=1}^{m}a_{j,\alpha}p_{j,\alpha} (0)\right)\varphi_1 \to 0,\quad \alpha\nearrow,
$$
because $\sum_{j=1}^{m}a_{j,\alpha}p_{j,\alpha} (0)=\Psi_{\alpha} (0)\to 0.$
It implies that
$$\Psi =\lim \Phi_{\alpha}, \quad \alpha\nearrow,
$$
 where $\Phi_{\alpha}=\Psi_{\alpha}-\psi_{\alpha} .$

Further, $\Phi_{\alpha} (0)=0$ and
\begin{multline*}
\frac{\Phi_{\alpha}}{z} =
\sum_{j=1}^{m} a_{j,\alpha}e^{-\mathrm{i}t_{j,\alpha}z} q_{j,\alpha}\varphi_j
+\sum_{j=1}^{m} a_{j,\alpha}p_{j,\alpha} (0)\frac{e^{-\mathrm{i}t_{j,\alpha}z}-1}{z} \varphi_j+\\
+\sum_{j=1}^{m} a_{j,\alpha}p_{j,\alpha} (0) \frac{\varphi_j-\varphi_1}{z} =\Sigma_1 +\Sigma_2+\Sigma_3,
\end{multline*}
where $q_{j,\alpha} =(p_{j,\alpha}-p_{j,\alpha} (0))z^{-1}$ are polynomials.
Clearly, $\Sigma_1\in\mathcal J$.

Let
 $\mathcal J_0$ be a submodule generated by 
$\varphi_1,$ $\dots ,$ $\varphi_m.$
Because of Corollary  \ref{cor-3}, $\mathcal J_0$ is  stable. 
 It follows that $\frac{\varphi_j-\varphi_1}{z}\in\mathcal J_0,$ $j=1,\dots , m.$  
Hence, $\Sigma_3\in \mathcal J_0\subset\mathcal J. $

Applying the results of \cite[VI.H.1]{Koosis-1}, we approximate every function
$$\frac{e^{-\mathrm{i}t_{j,\alpha }z}-1}{z}\varphi_j,\ \ j=1,\dots , m,$$
by the elements of 
$$\mathrm{span}\, \{e^{-\mathrm{i}tz}\varphi,\  t\in I\},$$
 and get $\Sigma_2\in\mathcal J.$

 Finally,
$$
\frac{\Phi_{\alpha}}{z}\in\mathcal J,\quad\text{and}\quad
\frac{\Psi}{z}=\lim \frac{\Phi_{\alpha}}{z} \in\mathcal J.
$$
That is,  $\mathcal J$ is stable at the origin.

Q.E.D.

\medskip

\begin{corollary}
Let 
 $I\subset (-a;a)$ be a relatively closed  interval, and
$\{S_{\alpha}\}\subset\mathcal U_a'$ be such that 
$\varphi_{\alpha}=\mathcal F (S_{\alpha})\in\mathcal M_{a}$, $\forall \alpha .$

Then, $W=\bigcap\limits_{\alpha} \mathrm{ker} \, T_{S_{\alpha}, I}$
 is a  weakly synthesable  subspace,  $I_W=I,$ and the sequence of exponents of
the system  $\mathrm{Exp} (W)$ equals
$\bigcap\limits_{\alpha}\mathcal Z_{\varphi_{\alpha}}.$
\label{cor-8}
\end{corollary}

Corollary \ref{cor-8} is justified by the same scheme as Corollary   \ref{cor-7} was done
(using Corollary \ref{cor-5} instead of Corollary \ref{cor-3}).

%%%%%%%%%%%%%%%%%%%%%%%%%%%%%%%%%%%%%%%%%%%%%%%%%%%%%%%%%%%%%%%%%%%%%%%%%%%%%%%%%%%%%%%%%%%%%%%%%%%%%%%%%%%%%%%
%%%%%%%%%%%%%%%%%%%%%%%%%%%%%%%%%%%%%%%%%%ACKNOWLEDGEMENTS
%%%%%%%%%%%%%%%%%%%%%%%%%%%%%%%%%%%%%%%%%%%%%%%%%%%%%%%%%%%%%%%%%%%%%%%%%%%%%%%%%%%%%%%%%%%%%

\smallskip

%\begin{acknowledgments}
{\sl 
The work is supported by Russian Science Foundation (code of scientific theme № 22-21-00026), https://rscf.ru/en/project/22-21-00026/.}
%\end{acknowledgments}

\def\bibname{\vspace*{-30mm}{

\vskip 1 cm \footnotesize
\begin{flushleft}
Natalia~Abuzyarova \\ % name of the authors from the same institution
Department of Mathematics and IT\\ % name of the department, where author works
Bashkir State University\\ % name of the university, where author works
32 Zaki Validi St,\\ % buseness address 
450076 Ufa, Russia; \\ % buseness address 
Institute of Mathematics with CC - Subdivision of the UFRC of RAS\\
 112 Chernyshevsky str. \\
  450008 Ufa, Russia;

E-mail: abnatf@gmail.com
\end{flushleft}

\end{document}